\documentclass[a4paper,10pt,oneside]{article}
\usepackage{graphicx}
\usepackage{amsmath, enumerate, amsfonts,amssymb, pictexwd, amscd}
\usepackage{amsthm}

\usepackage[latin1]{inputenc}
\usepackage[all]{xy}
\xyoption{arc}
\newtheorem{Conjecture}{Conjecture}[section]

\newtheorem{Remark}[Conjecture]{Remark}
\newtheorem{Theorem}[Conjecture]{Theorem}

\newtheorem{Proposition}[Conjecture]{Proposition}

\newcommand{\pcom}{^\wedge_p}

\begin{document}
\title{Group Models for Fusion Systems}
\author{Nora Seeliger}
\date{}
\maketitle
\begin{abstract}
 We study group models for fusion systems and construct homology decompositions for the models of Robinson and Leary-Stancu type. 
\end{abstract}
\section{Introduction}
In the topological theory
of $p-$local finite groups introduced by Broto, Levi, and Oliver
(see e.g. \cite{BLO2})
one tries to approximate the classifying space
of a finite group via the $p-$local structure of the group
as encoded by its fusion system of $p$-subgroups,
at least up to $\mathbb{F}_p-$cohomology.
The Martino-Priddy Conjecture, stated in \cite{MP} and proved in \cite{O1}, \cite{O2}, reveals that the homotopy type of the $p-$completed classifying space of a finite group is determined by $p-$fusion.
The $p$-local structure of a finite group is a special instance of a
fusion system, but not every fusion system arises in this way (see e.g.
\cite{RV}).

As in the case of groups, $p-$local finite groups have classifying spaces (see e.g. \cite{BLO2}), but these
need not be the classifying space of a finite group loc. cit., not even up to
mod $\mathbb{F}_p$-cohomology.
In this article we
study groups realizing fusion systems 
including the possibility of
perhaps even simultanously at different primes.
We construct homology decompositions to compute the
cohomology of the classifying spaces of group theory
models of fusion systems introduced by
Robinson and Leary-Stancu and investigate 
their homotopy types.

This work is part of the author's doctoral thesis at Aberdeen Institute of Mathematics under the supervision of Prof Ran Levi. I would like to thank my thesis advisor and the professors David Benson, Bill Dwyer, Radha Kessar, Assaf Libman, Markus Linckelmann, Bob Oliver and Geoffrey Robinson for kind advice and enriching discussions on the topic and moreover the referee for the generous help on the improvement of the exposition of the material.
\section{Preliminaries}

We review the basic definitions of fusion systems and centric linking systems and establish our notations. Our main references are \cite{BLO2}, \cite{BLO4} and \cite{IntroMarkus}.
Let $S$ be a finite $p$-group. A \textbf{fusion system} \textbf{$\mathcal{F}$} on $S$ is a category whose objects are all the subgroups of $S$, and which satisfies the following two properties for all $P,Q\leq S$: The set $Hom_{\mathcal{F}}(P,Q)$ contains injective group homomorphisms and, amongst them, all morphisms induced by conjugation of elements in $S$ and each element is the composite of an isomorphism in $\mathcal{F}$ followed by an inclusion. Two subgroups $P,Q\leq S$ will be called \textbf{$\mathcal{F}-$conjugate} if they are isomorphic in $\mathcal{F}$.
Define the group $Out_{\mathcal{F}}(P)=Aut_{\mathcal{F}}(P)/Inn(P)$, for all $P\leq S$, where $Inn(P)$ is the group of inner automorphisms of $P$.
 A subgroup $P\leq S$ is \textbf{fully centralized} resp. \textbf{fully normalized} in $\mathcal{F}$ if $|C_S(P)|\geq |C_S(P')|$ resp. $|N_S(P)|\geq |N_S(P')|$ for each subgroup $P'\leq S$ which is $\mathcal{F}$-conjugate to $P$. Every finite group gives rise to a canonical fusion system (see e.g.
\cite{BLO2}). The fusion system
 $\mathcal{F}$ is called \textbf{saturated} if, for each subgroup $P\leq S$ which is fully normalized in $\mathcal{F}$, $P$ is fully centralized in $\mathcal{F}$ and $Aut_S(P)\in Syl_p(Aut_{\mathcal{F}}(P))$, and, moreover, if $P\leq S$ and $\phi\in Hom_{\mathcal{F}}(P,S)$ are such that $\phi (P)$ is fully centralized, 
then there is a lift $\overline{\phi}\in Hom_{\mathcal{F}}(N_{\phi},S)$ such that $\overline{\phi}|_P=\phi$, where
$N_{\phi}=\{g\in N_S(P)|\phi c_g \phi^{-1}\in Aut_S(\phi (P))\}$. A subgroup $P\leq S$ will be called \textbf{$\mathcal{F}-$centric} if $C_S(P')\leq P'$ for all $P'$ which are $\mathcal{F}-$conjugate to $P$. A subgroup $P \leq S$ is called \textbf{$\mathcal{F}$-radical} if $O_p(Aut_{\mathcal{F}}(P))=Aut_P(P)$, where $O_p(-)$ denotes the maximal normal $p$-subgroup, and a subgroup $P \leq S$ is called \textbf{$\mathcal{F}$-essential} if $P$ is $\mathcal{F}$-centric and $Aut_{\mathcal{F}}(P)/Aut_P (P)$
has a strongly $p$-embedded subgroup. 
  Denote by \textbf{$\mathcal{F}^c$ }the full subcategory of $\mathcal{F}$ with objects the $\mathcal{F}-$centric subgroups of $S$.
	
 The \textbf{centric linking system associated to the saturated fusion system $\mathcal{F}$} is the category $\mathcal{L}$ whose objects are the
$\mathcal{F}$-centric subgroups of $S$, together with a functor
$\pi :\mathcal{L}\longrightarrow\mathcal{F}^c$,
and a \textbf{"distinguished" monomorphism $\delta _P:P\rightarrow Aut_{\mathcal{L}}(P)$} for each of the $\mathcal{F}$-centric subgroups $P\leq S$ such that the following conditions are satisfied. The functor $\pi$ is the identity on objects and surjective on morphisms. For each pair of objects $P,Q\in\mathcal{L}$ the group $Z(P)$ acts freely on $Mor_{\mathcal{L}}(P,Q)$ by precomposition (upon identifying the group $Z(P)$ with $\delta _P (Z(P))\leq Aut_{\mathcal{L}}(P)$), and the functor $\pi$ induces a bijection of sets 
$Mor_{\mathcal{L}}(P,Q)/Z(P)\overset{\simeq}{\longrightarrow}Hom_{\mathcal{F}}(P,Q).$
We have that for each of the $\mathcal{F}$-centric subgroups $P\leq S$ and each $x\in P$, $\pi (\delta _P(x))=c_x \in Aut_{\mathcal{F}}(P)$.
For every morphism $f\in Mor_{\mathcal{L}}(P,Q)$ and every element $x\in P$, $f\circ\delta_P(x)=\delta _Q(\pi f(x))\circ f$.

Let $\mathcal{F}$ and $\mathcal{F'}$ be fusion systems on finite $p$-groups $S$ and $S'$, respectively. A \texttt{}\textbf{{morphism of fusion systems }from $\mathcal{F}$ to $\mathcal{F'}$} is a pair $(\alpha, \Phi)$ consisting of a group homomophism $\alpha: S\rightarrow S'$, and a covariant functor $\Phi: \mathcal{F}\rightarrow\mathcal{F'}$ with the following properties: for any subgroup $Q$ of $S$, we have $\alpha (Q)=\Phi (Q)$ and for any morphism $\phi: Q\rightarrow R$ in $\mathcal{F}$, we have $\Phi(\phi)\circ\alpha |_Q=\alpha |_R\circ\phi$.

Let $G$ be a (possibly infinite) group. A finite subgroup $S$ of $G$ will be called a \textbf{Sylow $p-$subgroup} of $G$ if $S$ is a $p-$subgroup of $G$ and all $p-$subgroups of $G$ are conjugate to a subgroup of $S$.
A group $G$ is called \textbf{$p-$perfect} if $H_1(BG;\mathbb{Z}_p)=0$. Equivalently a group is $p-$perfect if it has no normal subgroup of index $p$. With this definition it is easy to see that a group generated by $p-$perfect subgroups is itself $p-$perfect. In particular, a group generated by $p'-$elements is $p-$perfect.  

Let $G_1$, and $G_2$ be groups with Sylow $p-$subgroups $S_1$, and $S_2$ respectively, and let $\phi :G_1\rightarrow G_2$ be a group homomorphism such that $\phi (S_1)\leq S_2$. The morphism $\phi$ will be called fusion preserving if $\phi |_{S_1}$ induces an isomorphism of fusion systems $\mathcal{F}_{S_1}(G_1)\cong\mathcal{F}_{S_2}(G_2)$.
Let $S$ be a finite $p-$group, and let $P_1,...,P_r,Q_1,...,Q_r$ be subgroups of $S$. Let $\phi _1,...,\phi _r$ be injective group homomorphisms $\phi _i:P_i\rightarrow Q_i$ for all $i=1,\cdots ,r$. The \textbf{fusion system generated by $\phi _1,...,\phi_r$} is the minimal fusion system $\mathcal{F}$ over $S$ containing $\phi _1,...,\phi _r$. 

Fix any pair $S\leq G$, where $G$ is a (possibly infinite) group and $S$ is a finite $p-$subgroup.
	Define $\mathcal{F}_S(G)$ to be the category whose objects are the subgroups of $S$, and where the set of morphisms between two subgroups $P$ and $Q$ fulfills $Mor_{\mathcal{F}_S(G)}(P,Q)=Hom_G(P,Q)=\{c_g\in Hom(P,Q)|g\in G, gPg^{-1}\leq Q\},$ where $c_g$ denotes the homomorphism conjugation by $g$ $(x\mapsto gxg^{-1})$. Note that $Hom_G(P,Q)\cong N_G(P,Q)/C_G(P),$ where $N_G(P,Q)=\{g\in G|gPg^{-1}\leq Q\}$ denotes the \textbf{\textbf{transporter set}}.
For each $P\leq S$, let $C'_G(P)$ be the maximal $p-$perfect subgroup of $C_G(P)$. Let $\mathcal{L}^c_S(G)$ be the category whose objects are the $\mathcal{F}_S(G)-$centric subgroups of $S$, and where $Mor_{\mathcal{L}^c_S(G)}(P,Q)=N_G(P,Q)/C_G'(P).$
	Let $\pi:\mathcal{L}^c_S(G)\rightarrow\mathcal{F}_S(G)$ be the functor which is the inclusion on objects and sends the class of $g\in N_G(P,Q)$ to conjugation by $g$. For each $\mathcal{F}_S(G)-$centric subgroup $P\leq G$, let $\delta _P:P\rightarrow Aut_{\mathcal{L}^c_S(G)}(P)$ be the monomorphism induced by the inclusion $P\leq N_G(P)$.

A triple $(S,\mathcal{F},\mathcal{L})$ where $S$ is a finite $p-$group, $\mathcal{F}$ is a saturated fusion system on $S$, and $\mathcal{L}$ is an associated centric linking system to $\mathcal{F}$, is called a \textbf{$p-$local finite group}. Its \textbf{classifying space} is $|\mathcal{L}|\pcom$ where $(-)\pcom$ denotes the $p-$completion functor in the sense of Bousfield and Kan, see \cite{BK}.
This is partly motivated by the fact that every finite group $G$ gives canonically rise to a $p-$local finite group $(S,\mathcal{F}_S(G),\mathcal{L}^c_S(G))$ and $BG\pcom\simeq |\mathcal{L}|\pcom$, see \cite{BLO2}. In particular, all fusion systems coming from finite groups are saturated.

Let $\mathcal{F}$ be a fusion system on the the finite $p-$group $S$. The fusion system $\mathcal{F}$ is called an\textbf{ Alperin fusion system} if there are subgroups $P_1=S,P_2,\cdots , P_r$ of $S$ and finite groups $L_1,\cdots ,L_r$ such that: $P_i$ is the largest normal $p-$subgroup of $L_i$, $C_{L_i}(P_i)=Z(P_i)$, and $L_i/P_i\cong Out_{\mathcal{F}}(P_i)$ for each $i$. Moreover $N_S(P_i)$ is a Sylow $p-$subgroup of $L_i$, and $\mathcal{F}_{N_S(P_i)}(L_i)$ is contained in $\mathcal{F}$ for each $i$ such that $\mathcal{F}$ is generated by all the $\mathcal{F}_{N_S(P_i)}(L_i)$. The groups $L_i$ are isomorphic to $Aut_{\mathcal{L}}(P_i)$ for all $i=1,\cdots ,r$. Recall that
every saturated fusion system is Alperin, as proven in \cite[Section 4]{BCGLO1}.

One can also define fusion systems and centric linking systems in a topological setting (see e. g. \cite[Definition 1.6]{BLO4}). We will need this when we make use of the fact that a group realizes a given fusion system if and only if its classifying space has a certain homotopy type. In particular we have for a $p-$local finite group $(S,\mathcal{F},\mathcal{L})$ and a group $G$ such that $\mathcal{F}_S(G)=\mathcal{F}$ that there is a map from the one-skeleton of the nerve of $\mathcal{L}$ to the classifying space : $|\mathcal{L}|^{(1)}\rightarrow BG$. 

Given a saturated fusion system $\mathcal{F}$ on a finite $p-$group $S$ it is not always true that there exists a finite group $G$, having $S$ as a Sylow $p-$subgroup, such that $\mathcal{F}_S=\mathcal{F}_S(G)$, (see \cite{BLO2}, chapter 9 for example). 
However for every fusion system $\mathcal{F}$ there exists an infinite group $\mathcal{G}$ such that $\mathcal{F}_S(\mathcal{G})=\mathcal{F}$. We now describe the constructions by G. Robinson \cite{Robinson1}, and I. Leary and R. Stancu \cite{Ian+Radu}.

Robinson's construction is specific to Alperin fusion systems. It is an iterated amalgam of the groups $L_i$ over the subgroups $N_S(P_i)$, where the groups $L_i$ and $P_i$, $i=1,\cdots ,r$ are appearing in the definition of an Alperin fusion system.

\begin{Theorem}[\cite{Robinson1}, Theorem 2.]
Let $\mathcal{F}$ be an Alperin fusion sytem on a finite $p-$group $S$ and associated groups $L_1,...,L_n$ as in the definiton. Then there exists a finitely generated group $\mathcal{G}$ which has $S$ as a Sylow $p-$subgroup such that the fusion system $\mathcal{F}$ is realized by $\mathcal{G}$. 
The group $\mathcal{G}$ is given explicitely by $\mathcal{G}=L_1\underset{N_S(P_2)}{*}L_2\underset{N_S(P_3)}{*}...\underset{N_S(P_n)}*L_n$ with $L_i,P_i$ as in the definiton.
\end{Theorem}
Corresponding to the various versions of Alperin's fusion theorem on a saturated fusion system $\mathcal{F}$ (reduction to automorphisms of $\mathcal{F}$-essential subgroups, of $\mathcal{F}$-centric subgroups respecetively, of $\mathcal{F}$-centric $\mathcal{F}$-radical subgroups) there are several choices for the groups generating $\mathcal{F}$ as an Alperin fusion system.

The infinite groups realizing arbitrary fusion systems constructed by I. Leary and R. Stancu are iterated HNN-constructions.
\begin{Theorem}[\cite{Ian+Radu}, Theorem 2.]
Suppose that $\mathcal{F}$ is the fusion system on $S$ generated by $\Phi=\{\phi_1, \cdots, \phi_r\}$. Let $T$ be a free group with free generators $t_1, \ldots, t_r$, and define $G$ as the quotient of the free product $S*T$ by the relations $t_i^{-1}ut_i=\phi_i(u)$ for all $i$ and for all $u\in P_i$. Then $S$ embeds as a $p-$Sylow subgroup of $G$ and $\mathcal{F}_S(G)=\mathcal{F}$. 
\end{Theorem}

\section{Group Models for Fusion Systems}
In this section we give a new construction of a group realizing a given fusion system. With it we give a construction of a group which realizes an arbitrary collection of fusion systems at different primes and study several examples.
\subsection{A functor to GROUP}
Let $p$ be a prime.
Define the category $FUSION(p)$. The objects of this category are fusion systems over finite $p-$groups and on morphisms we have morphisms between the respective fusion systems. Let $GROUP_{Syl_p}$ be the full subcategory of the category of groups where the objects are groups which have a Sylow $p-$subgroups. Define the functor $\textbf{F}: FUSION(p)\rightarrow GROUP_{Syl_p}$, as constructed in \cite[Corollary 4]{Ian+Radu}. Let $\mathcal{F}$ be an object of $FUSION(p)$, i. e. a fusion system over a finite $p-$group $S$. Then the functor $\textbf{F}$ takes $\mathcal{F}$ to the group $S*F(Mor(\mathcal{F}))/<\phi u\phi ^{-1}=\phi (u)\text{ }\forall \phi\in Mor(\mathcal{F}),\text{ }\phi : P\rightarrow Q,\text{ }u\in P>$, where $Mor(\mathcal{F})$ is the set of all morphisms in $\mathcal{F}$ and $F(Mor(\mathcal{F}))$ is the free group on the morphism set $Mor(\mathcal{F})$. Let $\mathcal{F}$ and $\mathcal{F'}$ be fusion systems over the finite $p-$groups $S$ and $S'$ respectively. Let  $(\alpha , \Phi ):\mathcal{F}\rightarrow\mathcal{F'}$ be a morphism of fusion systems between them. Define $\textbf{F}((\alpha, \Phi)):\textbf{F}(\mathcal{F})\longrightarrow\textbf{F}(\mathcal{F'})$ by
$s\mapsto\alpha (s)$ and $\phi\mapsto\Phi (\phi)$. The following remark is the same as \cite[Corollary 4]{Ian+Radu}.
\begin{Remark}The functor \textbf{F} is a left inverse to the canonical functor.
\end{Remark}
We give a new remark about the canonical functor in fusion theory.
\begin{Remark}
The canonical functor does not have a left adjoint.
\end{Remark} 
\underline{Proof}: Assume it had one. Then the functor which associates to a given group with a Sylow $p-$subgroup its canonical fusion system will be a right adjoint and therefore preserve pullback diagrams. Recall that there are two nonisomorphic fusion systems over the cyclic group of order $3$ $C_3$: the trivial one and $\mathcal{F}_{C_3}(\Sigma _3)$. The fusion system of the pullback $\Sigma _3\rightarrow\Sigma _3\underset{C_3}{*}\Sigma _3\leftarrow\Sigma _3$ is not isomorphic to the pullback $\mathcal{F}_{C_3}(\Sigma _3)\rightarrow\mathcal{F}_{C_3}(\Sigma _3\underset{C_3}{*}\Sigma _3)\leftarrow\mathcal{F}_{C_3}(\Sigma _3)$. $\Box $\\[0.3cm]
Let $\mathcal{F}$ be a fusion system over a finite $p-$group $S$.
There is a group model $\mathcal{G}$ such that $\mathcal{F}_S(\mathcal{G})=\mathcal{F}$ which comes with the property that for every group $G$ such that $\mathcal{F}_S(G)=\mathcal{F}$ there exists $\mathcal{G}$ of this type such that $\mathcal{G}$ surjects on $G$ in a fusion preserving way. This result is \cite[Corollary 3]{Ian+Radu}.

\subsection{Graphs of Groups}
We construct a group realizing an arbitrary collection of fusion systems at different primes. To begin, a short resume from graphs of groups we need.
A \textbf{finite directed graph} $\Gamma$ consists of two sets, the \textbf{vertices} $V$ and the \textbf{directed edges} $E$, together with two functions $\iota, \tau : E\rightarrow V$. For $e \in E,\iota (e)$ is called the \textbf{initial vertex} of $e$ and $\tau (e)$ is the \textbf{terminal vertex} of $e$. Multiple edges and loops are allowed in this definition. The graph $\Gamma$ is \textbf{connected} if the only equivalence relation on $V$ that contains all pairs $(\iota (e),\tau (e))$ is the relation with just one class. 
A graph $\Gamma$ may be viewed as a category, with objects the disjoint union of $V$ and $E$ and two non-identity morphisms with domain $e$ for each $e\in E$, one morphism $e\rightarrow \iota (e)$ and one morphism $e\rightarrow \tau (e)$.
A \textbf{graph $\Gamma$ of groups} is a connected graph $\Gamma$ together with groups $G_v, G_e$ for each vertex and edge and injective group homomorphims $f_{e,\iota}:G_e\rightarrow G_{\iota (e)}$ and $f_{e,\tau (e)}:G_e\rightarrow G_{\tau (e)}$ for each edge $e$. If a graph is seen as a category, then a graph of groups can be seen as a functor from that category to the category of groups with injective homomorphisms.
The following theorem is of particular interest since the exoticity of the Solomon fusion systems \cite{Solomon} was shown (in today's language), via the noncompatibility of fusion systems at different primes in a finite group.
\begin{Theorem}
Let $p_1,\cdots p_m$ be a collection of different primes, $S_1,\cdots, S_m$ be a collection of $p_i-$groups respectively and $\mathcal{F}_{S_i}$ a collection of fusion systems over $S_i$ respectively. Then there exists a group $G$ such that $S_i\in Syl_{p_i}(G)$ and $\mathcal{F}_{S_i}(G)=\mathcal{F}_{S_i}$ for all $i$.
\end{Theorem}
\underline{Proof:}  
Let $\mathcal{G}_i$ be models of Leary-Stancu-type for the $\mathcal{F}_{S_i}$ for all $i$ respectively. Define  $G:=\overset{m}{\underset{i=1}{\times}}\mathcal{G}_i$. $S_i\in Syl_{p_i}(G)$ since all finite subgroups of $\mathcal{G}_i$ are conjugate to $S_i$ for all $i$ and therefore all finite $p_i-$groups are conjugate to a subgroup of $S_i$ for all $i$. Obviously $\mathcal{F}_{S_i}\subseteq \mathcal{F}_{S_i}(G)$. Since all elements of $\mathcal{G}_j$, $i\neq j$, commute with $\mathcal{G}_i$ we obtain $\mathcal{F}_{S_i}= \mathcal{F}_{S_i}(G)$ for all $i$. More generally we have the following. Let $\mathcal{G}$ be a group having a Sylow $p-$subgroup $S$ and $\mathcal{H}$ be a group that does not contain any $p-$element. Then $S$ is a Sylow $p-$subgroup of $\mathcal{G}\times\mathcal{H}$ and $\mathcal{F}_S(\mathcal{G})=\mathcal{F}_S(\mathcal{G}\times\mathcal{H})$. The same result is true if we replace $\mathcal{G}\times\mathcal{H}$ by $\mathcal{G}*\mathcal{H}$. This last statement is a consequence of \cite[Theorem 5]{Ian+Radu}.$\Box$
\begin{Remark}
Instead of taking the direct product of the $\mathcal{G}_i$ it is also possible to take the free product $G$ of the $\mathcal{G}_i$. We have $S_i\in Syl_{p_i}(G)$ and $\mathcal{F}_{S_i}\subseteq \mathcal{F}_{S_i}(G)$ for all $i$ as before. Since the only element of $S_i$ that is conjugate into $S_i$ by $G_j$, $i\neq j$, is the trivial. Therefore we obtain again $\mathcal{F}_{S_i}= \mathcal{F}_{S_i}(G)$ for all $i$. Note that direct product and amalgamated product of respective Robinson models will do as well as long as the orders of $L_1,...,L_n$ of each factor are coprime to all the remaining primes $p_i$. 
\end{Remark}
The following remark is a particular case of \cite[Theorem 5]{Ian+Radu}.
\begin{Remark} Let $\mathcal{F}$ be a fusion system over the finite $p-$group $S$.
Let $G, G'$ be groups such that $S \in Syl_p(G),S \in Syl_p(G'),\mathcal{F}_S(G)=\mathcal{F}$, and $\mathcal{F}_S(G')=\mathcal{F}$. Let $\mathcal{G}:=G\underset{S}{*}G'$. Then $S\in Syl_p(\mathcal{G})$ and $\mathcal{F}=\mathcal{F}_S(\mathcal{G})$.
\end{Remark}
\section{Group Models and Homology Decompositions}  
In \cite{CE} Cartan and Eilenberg show that the $\mathbb{F}_p-$cohomology ring of a finite group is given as the subring of stable elements of the cohomology ring of the Sylow $p-$subgoup. In \cite{BLO2} Broto, Levi and Oliver extend the result to the cohomology of the classifying space of the $p-$local finite group. We now relate their results to the cohomology of the group models studied by Robinson and Leary-Stancu. Before we do this we give a general result relating the cohomology of a group model and the ring of stable elements. Recall that the ring of $\mathcal{F}-$stable elements $H^*(\mathcal{F})$ is isomorphic to $H^*(|\mathcal{L}|)$ by \cite[Theorem 5.8]{BLO2}.
\begin{Theorem}
Let $(S,\mathcal{F},\mathcal{L})$ be a $p-$local finite group and $
\mathcal{G}$
a discrete group such that $S\in Syl_p(\mathcal{G})$ and $\mathcal{F}=\mathcal{F}_S(\mathcal{G})$ . Then
there exist a natural map of algebras $H^*(B\mathcal{G})\overset{q}{\rightarrow} H^*(\mathcal{F})$ making $H^*(\mathcal{F})$ a module over $H^*(B\mathcal{G})$.
\end{Theorem}
\underline{Proof:} We show that the restriction map $Res^{\mathcal{G}}_S:H^*(B\mathcal{G})\rightarrow H^*(BS)$ factors through the ring of stable elements $H^*(\mathcal{F})\subset H^*(BS)$. For every subgroup $P$ of $S$ we have a map $Res^{\mathcal{G}}_P: H^*(B\mathcal{G})\rightarrow H^*(BP)$. 
Since we have $\mathcal{F}=\mathcal{F}_S(\mathcal{G})$ we obtain that for all subgroups $P,  Q\leq S$ and all morphisms $\phi\in Mor_{\mathcal{F}}(P,Q)$ that $\phi =c_g$ for some $g\in \mathcal{G}$. Therefore the outer triangle commutes 
\xymatrix@R=40pt@C=40pt{ &&{H^*(BP)}&{}\\
{H^*(B\mathcal{G})}\ar[1,2]_{Res^{\mathcal{G}}_Q}\ar[-1,2]^{Res^{\mathcal{G}}_P}\ar[0,1]^{Res^{\mathcal{G}}_S}&{H^*(BS)}\ar[1,1]^{Res^S_Q}\ar[-1,1]_{Res^S_P}&&\\
&&{H^*(BQ)}\ar[-2,0]^{\phi ^*}\restore&{}\\
} \\
since the diagram
\xymatrix@R=7pt@C=7pt{ {BP}\ar[0,2]^{B\phi}\ar[1,1]_{Bincl}&{}&{BQ}\ar[1,-1]^{Bincl}\\
&{B\mathcal{G}}&\\} 
commutes up to homotopy and the map $Res^{\mathcal{G}}_S:H^*(B\mathcal{G})\rightarrow H^*(BS)$ factors through the ring of stable elements giving rise to the map $q:H^*(B\mathcal{G})\rightarrow H^*(\mathcal{F})$.$\Box$
\subsection{Homology Decomposition for Robinson's Models} 
We investigate the cohomology of Robinson's models for a saturated fusion system $\mathcal{F}$ over a finite $p$-group $S$. Let $P_1,...,P_n$ be $\mathcal{F}-$centric subgroups, fully normalized in $\mathcal{F}$ such that $P_1=S$ and $\mathcal{F}$ is generated by $Aut_{\mathcal{F}}(P_1), ..., Aut_{\mathcal{F}}(P_n)$. For each $P_i$, let $L_i=Aut_{\mathcal{L}}(P_i)$ be as in the definition of Alperin fusion systems. We recall that such $L_i`s$ exist and are unique (see e.g. \cite[Section 4]{BCGLO1}), and moreover $L_i=Aut_{\mathcal{L}}(P_i)$. In the following $\mathcal{G}$ will always be a model for $\mathcal{F}$ of Robinson type, i. e.
$\mathcal{G}=L_1\underset{N_S(P_2)}{*}L_2\underset{N_S(P_3)}{*}...\underset{N_S(P_n)}*L_n$.
\begin{Theorem}
Let $(S,\mathcal{F},\mathcal{L})$ be a $p-$local finite group and 
$\mathcal{G}$
the model of Robinson type for $\mathcal{F}$ presented above. Then
there exist natural maps of algebras over the Steenrod algebra $H^*(B\mathcal{G})\overset{q}{\rightarrow} H^*(|\mathcal{L}|)$ and $H^*(|\mathcal{L}|)\overset{r^*}{\rightarrow} H^*(B\mathcal{G})$ such that we obtain a split short exact sequence of unstable modules over the Steenrod algebra\begin{eqnarray*}0\longrightarrow W\overset{incl}{\longrightarrow }H^*(B\mathcal{G})\overset{\overset{r^*}{\leftarrow}}{\underset{q}{\longrightarrow}}H^*(|\mathcal{L}|)\longrightarrow 0,\end{eqnarray*} where $W\cong Ker(Res^{\mathcal{G}}_S)$. 
\end{Theorem}
\underline{Proof:}
Let $\mathcal{C}$ be the following category. \\
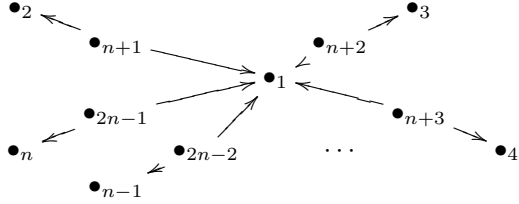
\xymatrix@R=1pt@C=3pt{
{\bullet _2}&&&&&&{\bullet _3}&&\\ 
&&{\bullet _{n+1}}\ar[1,2]\ar[-1,-2]&&&{\bullet _{n+2}}\ar[1,-1]\ar[-1,1]&&&\\
&&&&{\bullet _1}&&&&\\
&&{\bullet _{2n-1}}\ar[1,-2]\ar[-1,2]&&&&{\bullet _{n+3}}\ar[-1,-2]\ar[1,2]&&&\\
{\bullet _n}&&&{\bullet _{2n-2}}\ar[-2,1]\ar[1,-1]&&{\cdots}&&&{\bullet _4}&\\
&&{\bullet _{n-1}}&&&&&&&\\
} \\
Denote by $\phi _{i,j}:\bullet_{i}\rightarrow \bullet _j$ the unique morphism in $\mathcal{C}$ between $\bullet _i$ and $\bullet _j$ if it exists.
Let $F:\mathcal{C}\rightarrow Spaces$ be a functor with $F(\bullet _i)=BL_i$ for $i=1,...,n,  BN_S(P_i)$ for $i=n+1,...,2n-1$ and $F(\phi _{i,j})=Bincl : F(\bullet _i)\rightarrow F(\bullet _j)$ for all $\phi _{i,j}:\bullet _i\rightarrow \bullet _j$ in $\mathcal{C}$, $i=n+1,...,2n-1$, $j=i-n+1,1$. Note that $\underset{\mathcal{C}}{hocolim(F)}$ is a $K(\mathcal{G},1)$ and the graph of groups we use here is homotopy equivalent to the star-shaped graph of groups constructed in the proof of \cite[Theorem 5]{Ian+Radu}. Since $L_i=Aut_{\mathcal{L}}(P_i)$ for all $i=1,...,n$ we have a functor $\mathcal{B}L_i\rightarrow\mathcal{L}$ which sends the unique object $\bullet$ to $P_i$ and a morphism $x$ to the corresponding morphism in $Aut_{\mathcal{L}}(P_i)$ for all $i=1,...,n$. Therefore we obtain a map $BL_i$ to $|\mathcal{L}|$ for all $i=1,...,n$. Note that all the diagrams
\xymatrix@R=3pt@C=3pt{
&{BL_1}\ar[1,1]^{}&{}\\
{BN_S(P_i})\ar[-1,1]^{Bincl}\ar[1,1]_{Bincl}&&{|\mathcal{L}|}\\
&{BL_i}\ar[-1,1]_{}\restore&{}
}
commute up to homotopy since the existence of the linking system guarantees that we can find a compatible system of lifts of the inclusion $\iota _{N_S(P_i),S}$ in $\mathcal{L}$ for all $i=1,...,n$ such that all the diagrams
\xymatrix@R=3pt@C=3pt{ 
&{\mathcal{B}L_1}\ar[1,1]^{}&{}\\
{\mathcal{B}N_S(P_i})\ar[-1,1]^{\mathcal{B}incl}\ar[1,1]_{\mathcal{B}incl}&&{\mathcal{L}}\\
&{\mathcal{B}L_i}\ar[-1,1]_{}\restore&{}\\
}
commute up to the natural transformation
which takes the object $\bullet\in Obj(\mathcal{B}N_S(P_i))$ to $\iota _{N_S(P_i),S}$ for $i=1,...n.$
We obtain a map from the 1-skeleton of the homotopy colimit of the functor $F$ over the category $\mathcal{C}$ to $|\mathcal{L}|$. Since $\mathcal{C}$ is a $1-$dimensional category we obtain a map from $B\mathcal{G}$ to $|\mathcal{L}|$. This map will be denoted by $r$ inducing $H^*(|\mathcal{L}|)\overset{r^*}{\rightarrow } H^*(B\mathcal{G})$.\\
From the construction of the isomorphism $H^*(|\mathcal{L}|)\rightarrow H^*(\mathcal{F})$ in \cite[Theorem 5.8]{BLO2} it follows that the natural map $BS\stackrel{B(\delta _S)}{\longrightarrow } |\mathcal{L}|$ induces an isomorphism in cohomology on the image $H^*(|\mathcal{L}|)\stackrel{B(\delta _S)^*}{\longrightarrow } H^*(\mathcal{F})\subset H^*(BS). $ We have the following commutative diagram
\xymatrix@R=7pt@C=7pt{ {BS}\ar[0,2]^{B(\delta _S)}\ar[1,1]_{Bincl}&{}&{|\mathcal{L}|}\\
&{B\mathcal{G}}\ar[-1,1]_{r}&\\} inducing a commutative diagram of unstable algebras over the Steenrod algebra
\xymatrix@R=7pt@C=7pt{ {H^*(BS)}&{}&{H^*(|\mathcal{L}|)}\ar[0,-2]_{B(\delta _S)^*}\ar[1,-1]^{r^*}\\
&{H^*(B\mathcal{G})}\ar[-1,-1]^{q}&\\} which shows that $q\circ r^*=id_{H^*(|\mathcal{L}|)}$ where $H^*(|\mathcal{L}|)$ is identified with $H^*(\mathcal{F})$ via the isomorphism quoted above.

Denote the kernel of the map $q$ by $W$. We have the following commutative diagram of unstable algebras 
\xymatrix@R=9pt@C=9pt{ {H^*(|\mathcal{L}|)}\ar[1,1]_{incl}\ar[0,1]^{r^*}_{\underset{q}{\longleftarrow}}&{H^*(B\mathcal{G})}\ar[1,0]^{Res^{\mathcal{G}}_S}\\
&{H^*(BS).}\\
}Commutativity implies that $W\cong Ker(Res^{\mathcal{G}}_S)$ in the category of unstable modules.$\Box$

\begin{Theorem}
Let $\mathcal{F}$ be an Alperin fusion system and $\mathcal{G}$ a model of Robinson type for it. Then $B\mathcal{G}$ is $p-$good.
\end{Theorem}
\underline{Proof:} 
The group $\mathcal{G}$ is a finite amalgam of finite groups. Note that each $L_i$ is generated by $N_S(P_i)$ and elements of $p'-$order. Therefore $\mathcal{G}$ is generated by elements of $p'-$order and $S$. Let $K$ be the subgroup of $\mathcal{G}$ generated by all elements of $p'-$order. Note that $K\triangleleft\mathcal{G}$ and $S$ surjects on $\mathcal{G}/K$ and therefore $\mathcal{G}/K$ is a finite $p-$group. The group $K$ is $p-$perfect since it is generated by $p'-$elements and therefore we have $H^1(BK;\mathbb{F}_p)=0$. Let $X$ be the cover of $B\mathcal{G}$ with fundamental group $K$. Then, using \cite[VII.3.2]{BK}, we have that $X$ is $p-$good and $X\pcom$ is simply connected. Hence $X\pcom\rightarrow B\mathcal{G}\pcom\rightarrow B(\mathcal{G}/K)$ is a fibration sequence and $B\mathcal{G}\pcom$ is $p-$complete by \cite[II.5.2(iv)]{BK}. So $B\mathcal{G}$ is $p-$good. $\Box$
\begin{Theorem}
Let $(S, \mathcal{F}, \mathcal{L})$ be a $p-$local finite group and $\mathcal{G}$ be a model of Robinson type for it. Then $H^*(B\mathcal{G})$ is finitely generated.
\end{Theorem}
\underline{Proof:}
 Note that we have a map $B\mathcal{G}=\underset{\mathcal{C}}{hocolim(F)}\rightarrow |\mathcal{L}|$ where $F$ and $\mathcal{C}$ are as defined in the proof of Theorem 4.2. for the model of Robinson type $\mathcal{G}$. Note that $N_S(P_i)\in Syl_p(L_i)$ for all $i=1,...,n$. It follows from \cite[Lemma 2.3.]{BLO1} and \cite[Theorem 4.4.(a)]{BLO2} that $H^*(BP_i)$ is finitely generated over $H^*(|\mathcal{L}|)$ for all $i=1,...,n$ since $H^*(BP_i)$ is finitely generated over $H^*(BS\pcom)$ and therefore over $H^*(BS)$ and therefore over $H^*(|\mathcal{L}|\pcom)$ which is contained in $H^*(BS)$ as a subring and contains the image of the map $H^*(BP_i)\rightarrow H^*(BS)$. Moreover we have that $H^*(|\mathcal{L}|)$ is noetherian as follows from \cite[Proposition 1.1. and Theorem 5.8.]{BLO2}. Therefore the Bousfield-Kan spectral sequence for $H^*(B\mathcal{G})$ is a spectral sequence of finitely generated $H^*(|\mathcal{L}|)-$modules, the $E_2$ term with $E_2^{s,t}= \underset{\mathcal{C}}{lim^s}H^t(F(-);\mathbb{F}_p)$ is concentrated in the first two columns and $E_{2}=E_{\infty}$ for placement reasons, see \cite{Dwyer} for a reference for the Bousfield-Kan spectral sequence of a homotopy colimit. Therefore $H^*(B\mathcal{G})$ is a finitely generated module over $H^*(|\mathcal{L}|)$ and in particular noetherian. $\Box$ \\

A space $A$ is \textbf{a stable retract} of a space $X$ if there exists a map $f:A\rightarrow X$ such that the induced map on suspension spectra $\Sigma ^{\infty}f:\Sigma ^{\infty}A\rightarrow \Sigma ^{\infty}X$ is a retract.
\begin{Theorem}
Let $(S,\mathcal{F},\mathcal{L})$ be a $p-$local finite group and $\mathcal{G}$ a model of Robinson type for $\mathcal{F}$. Then
$|\mathcal{L}|\pcom$ is a stable retract of $B\mathcal{G}\pcom$.
\end{Theorem}
\underline{Proof:} 
The diagram
\xymatrix@R=9pt@C=9pt{&{\Sigma ^{\infty}BS\pcom}\ar[1,1]^{\Sigma ^{\infty}Bincl\pcom}\ar[1,-1]_{\Sigma ^{\infty}B(\delta _S)\pcom}&{}\\
{\Sigma ^{\infty}|\mathcal{L}|\pcom}&&{\Sigma ^{\infty}B\mathcal{G}\pcom}\ar[0,-2]^{\Sigma ^{\infty}r\pcom}\\
}
commutes where $r$ is the map constructed in the proof of Theorem $4.1$.
By Ragnarsson's work \cite{Ragnarsson} there is a map $\sigma _{\mathcal{F}}:\Sigma ^{\infty}|\mathcal{L}|\pcom\rightarrow\Sigma ^{\infty}BS \pcom$ such that the composition of maps $ \Sigma ^{\infty}|\mathcal{L}|\pcom\overset{\sigma _{\mathcal{F}}}{\longrightarrow}\Sigma ^{\infty}BS\pcom\overset{\Sigma ^{\infty}B(\delta _S)\pcom}{\longrightarrow}\Sigma ^{\infty}|\mathcal{L}|\pcom$ is the identity.
Since $\Sigma ^{\infty}B(\delta _S)\pcom\circ \sigma _{\mathcal{F}} =\Sigma ^{\infty}r\pcom\circ\Sigma ^{\infty}Bincl\pcom\circ \sigma _{\mathcal{F}}$ we have $|\mathcal{L}|\pcom$ is a stable retract of $B\mathcal{G}\pcom$.$\Box$\\

We outline two examples of a fusion system and a model of Robinson type for it: In the first the cohomology of the model is isomorphic to the stable elements and in the second it is not.
\begin{Proposition}
Let $G=PSL_2(7)$ be the projective special linear group of rank $2$ over the field of $7$ elements.
There exists $\mathcal{G}$ a model of Robinson type for the $2-$local finite group associated to $G$ such that $H^*(B\mathcal{G};\mathbb{F}_2)\cong H^*(|\mathcal{L}|;\mathbb{F}_2)$.
\end{Proposition}
\underline{Proof:} Note that $PSL_2(7)$ has a $2-$Sylow subgroup $D_8$, isomorphic to the dihedral group of order $8$. Note that $\mathcal{F}=\mathcal{F}_{D_8}(PSL_2(7))$. Denote by $V$ and $W$ two representatives of the respective conjugacy classes of $\mathcal{F}-$centric $\mathcal{F}-$radicals of $D_8$. We know from \cite[Example 8.8]{IntroMarkus} that $Aut_{\mathcal{F}}(V)=Aut_{\mathcal{F}}(W)=\Sigma _3$ and therefore we have that $Aut_{\mathcal{L}}(V)=Aut_{\mathcal{L}}(W)=\Sigma _4$ and $\mathcal{G} =\Sigma _4\underset{D_8}{*}{\Sigma _ 4}$ is a model of Robinson type for the $2-$local finite group associated to $G$. It follows from \cite[Theorem 3.5]{BCGLO1} that  $B(D_8\underset{D_8}{*}\Sigma _4\underset{D_8}{*}\Sigma _4)$ is weakly equivalent to $|\mathcal{L}|$ before completion and therefore we obtain that $H^*(B\mathcal{G};\mathbb{F}_2)\cong H^*(|\mathcal{L}|;\mathbb{F}_2)$.$\Box$
\begin{Proposition}
Let $p$ be an odd prime. Let
$C_p \wr C_p\in Syl_p (\Sigma _{p^2})$
and $(S,\mathcal{F},\mathcal{L})$ be the associated $p-$local finite group. Then there does not exist a model of Robinson type for $\mathcal{F}$ such that $H^*(B\mathcal{G})$ is isomorphic to $H^*(|\mathcal{L}|)$.
\end{Proposition}
\underline{Proof:}
We know from \cite[Section 2]{AF} that representatives of isomorphism classes of $\mathcal{F}-$centric $\mathcal{F}-$radical subgroups of $C_p\wr C_p$ are
the Sylow $p-$subgroup $S=C_p\wr C_p$ where the normalizer in $\Sigma _{p^2}$ is $C_p\wr C_p\rtimes (GL_1(\mathbb{F}_p)\times GL_1(\mathbb{F}_p))$ with the diagonal action on $C_p$, $P_2=C_{p}\times C_p$ embedded via its action on itself via translation, the normalizer is $P_2\rtimes GL_2(\mathbb{F}_p)$, and $P_3=C_p^{\times p}$ as a subgroup of $\Sigma _p^{\times p}\leq\Sigma _{p^2}$, where the normalizer in $\Sigma _{p^2}$ is $(C_p\rtimes GL_1(\mathbb{F}_p))\wr\Sigma _p.$ 
Denote $L_i=Aut_{\mathcal{L}}(P_i)$ for $i=1,2,3$. The centralizers are contained in the centers.
The model of Robinson type we consider is $\mathcal{G}=L_1\underset{N_S(P_2)}{*}L_2\underset{N_S(P_3)}{*}L_{3}$.
We prove that $H^1(L_1;\mathbb{F}_p)=0$. Let $S:=C_p\wr C_p$ and $G:=S\rtimes (GL_1(\mathbb{F}_p)\times GL_1(\mathbb{F}_p))$ where $(GL_1(\mathbb{F}_p)\times GL_1(\mathbb{F}_p))$ acts diagonally on $V$. Assume $\phi :G\rightarrow C_p$ is a surjective group homomorphism. Let $K:=Ker\phi$. Then $|K|=p^p(p-1)^2$. We have $[S,K]\leq S\cap K$ since $S\triangleleft G$ and $K\triangleleft G$. Moreover $|S\cap K|=p^p$. The quotient $S/S\cap K$ is centralised via $K$ and therefore via $G$ and therefore via $GL_1(\mathbb{F}_p)\times GL_1(\mathbb{F}_p)$. This is a contradiction and $\phi$ cannot exist.
Analogously it follows that $H^1(L_2;\mathbb{F}_p)=0$.
We also have $H^1(L_3)=0$ since $L_3=(C_p\rtimes C_{p-1})\wr\Sigma _p$, and the abelianisation of $\Sigma _p$ is $C_2$ and since $p>2$ it follows from the Serre spectral sequence for the extension $0\rightarrow (C_p\rtimes C_{p-1})^p\rightarrow (C_p\rtimes C_{p-1})^p\rtimes C_2\rightarrow C_2\rightarrow 0$ that $H^1(L_3)=0$. 
We have $H^1(N_{P_1}(C_{p}\times C_p);\mathbb{F}_p)\neq 0$ and moreover we obtain $N_S(P_3)=S$.
In the corresponding Mayer-Vietoris-Sequence where the covering spaces are $BL_1$ and $BL_i$, for $i=2,...,n$ and the intersections are the corresponding $BN_S(P_i)$ we obtain\\
\xymatrix@R=3pt@C=4pt{ 
&&&{0}\ar@{=}[1,0]&{}&&&&&&&&\\
&&&{H^1(L_1)}&{}&&&&&{H^2(L_1)}&&&\\
&&&{\times}&{H^1(N_{S}(P_2))}&&&&&{\times}&&&\\
{H^1(B\mathcal{G})}\ar[0,3]&&&{H^1(L_2)}\ar[0,1]&{\times}\ar[0,2]&&{H^2(B\mathcal{G})}\ar[0,3]^{f}&&&{H^2(L_2)}\ar[0,3]&&&{...}\\
&&&{\times}&{H^1(N_{S}(P_3))}&&&&&{\times}&&&\\
&&&{H^1(L_3)}&&&&&&{H^2(L_3)}&&&\\
&&&{0}\ar@{=}[-1,0]&&&&&&&&&\\
}\\
Since $0\neq H^1(S)\times H^1(N_S(P_2))$ is the kernel of the map $f$ and $H^2(\Sigma _{p^2};\mathbb{F}_p)=0$ we conclude that the cohomology of $B\mathcal{G}$ is strictly bigger than the stable elements. Any other model $\mathcal{G}'$ of Robinson type contains $L_1$ and groups isomorphic to $L_2$ and $L_3$ so there cannot exist a model of Robinson type which has $\mathbb{F}_p-$cohomology isomorphic to the stable elements ruining the hope that we can always find a model of Robinson type for which the cohomology is isomorphic to the stable elements.$\Box$
\subsection{Homology Decompositions for Leary-Stancu groups} 
We also establish analogous results for
the group constructed by I. Leary and R. Stancu.
\begin{Theorem}
Let $(S,\mathcal{F},\mathcal{L})$ be a $p-$local finite group and 
$\mathcal{G}$
a model of Leary-Stancu type for $\mathcal{F}$. Then
there exist natural maps of algebras over the Steenrod algebra $H^*(B\mathcal{G})\overset{q}{\rightarrow} H^*(|\mathcal{L}|)$ and $H^*(|\mathcal{L}|)\overset{s^*}{\rightarrow} H^*(B\mathcal{G})$ such that we obtain a split short exact sequence of unstable modules over the Steenrod algebra\begin{eqnarray*}0\longrightarrow W\overset{incl}{\longrightarrow }H^*(B\mathcal{G})\overset{\overset{r^*}{\leftarrow}}{\underset{q}{\longrightarrow}}H^*(|\mathcal{L}|)\longrightarrow 0,\end{eqnarray*} where $W\cong Ker(Res^{\mathcal{G}}_S)$. 
\end{Theorem}
\underline{Proof:} Let $ \{\phi _1,..., \phi _n \} $ with $\phi _i :P_i\rightarrow Q_i$ for all $i=1,...,n$ be the set of morphisms used in the construction of $\mathcal{G}$. We will consider them from this point on as morphisms $\phi _i:P_i\rightarrow S$ for all $i=1,...,n$. In \cite{Ian+Radu} the authors show that $B\mathcal{G}=\underset{\mathcal{D}}{hocolim F}$ where $\mathcal{D}$ is the following category:\\
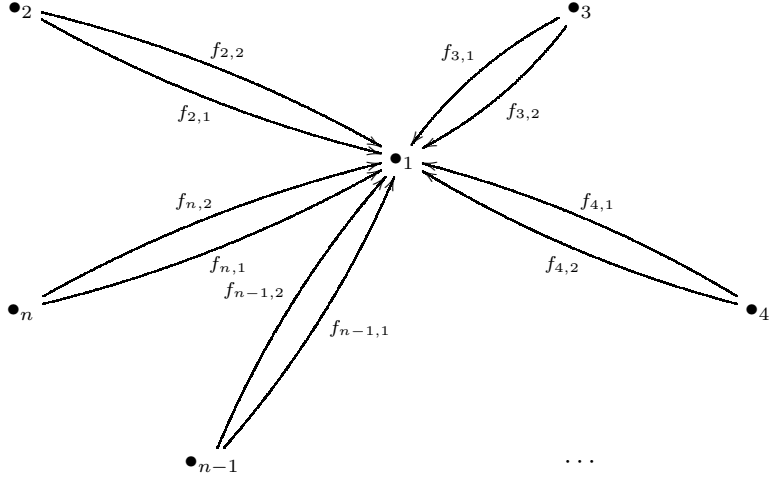
\xymatrix@R=45pt@C=50pt{
{\bullet _2}\ar @/^/[1,2]^{f_{2,2}}\ar @/_/[1,2]_{f_{2,1}}&&&{\bullet _3}\ar @/^/[1,-1]^{f_{3,2}}\ar @/_/[1,-1]_{f_{3,1}}&\\
&&{\bullet _1}&&\\
{\bullet _n}\ar @/_/[-1,2]_{f_{n,1}}\ar @/^/[-1,2]^{f_{n,2}}&&&&{\bullet _4}\ar @/^/[-1,-2]^{f_{4,2}}\ar @/_/[-1,-2]_{f_{4,1}}\\
&{\bullet _{n-1}}\ar @/_/[-2,1]_{f_{n-1,1}}\ar @/^/[-2,1]^{f_{n-1,2}}&&{\cdots}&\\
}\\
and $F$ is a functor to spaces with $F(\bullet _1)=BS$ and $F(\bullet _i)=BP_i$ and $F(f_{i1})=Bincl:BP_i\rightarrow BS$ and $F(f_{i2})=B\phi _i:BP_i\rightarrow BS$ for all $i=1,...,n$.\\
Due to Alperin's fusion theorem  there exist for all $i=1,\cdots, n$ an index $k(i)$ and a family of $\mathcal{F}-$centric subgoups $P^1_i,\cdots ,P^{k(i)}_i$ and for all $j=1,\cdots,k(i)$ $\psi ^j_i\in Aut_{\mathcal{F}}(P^j_i)$ such that for all $i=1,\cdots ,n$, for all $x\in P_i$ we have $\phi _i(x)=\psi ^{k(i)}_i\circ \psi ^{k(i)-1}_i\circ\cdots\circ\psi ^1_i(x)$. Whenever we consider the automorphisms $\psi ^j_i$ as morphisms $\psi ^j_i:P^j_i\rightarrow S$ they will be denoted $\widehat{\psi ^j_i}:P^j_i\rightarrow S$ for all $i=1,\cdots ,n$ and $j=1,\cdots ,k(i)$. Note that all the groups $P_i^j$ are $\mathcal{F}-$centric subgroups of $S$ and therefore there is a functor $F_i^{2j}:\mathcal{B}P^j_i\rightarrow\mathcal{L}=incl\circ\delta _{P_i^j}$ for all $i=1,...,n$, $j=1,\cdots ,k(i)$. Therefore we have a map $BP_i^j\rightarrow |\mathcal{L}|$ for all $i=1,...,n$, $j=1,\cdots ,k(i)$. 
Each object $BP$ of the image of the category $\mathcal{D}$ under $F$ gets mapped to $BAut_{\mathcal{L}}(P)\subset |\mathcal{L}|$. For all $i=1,\cdots ,n$ define the family of functors $\{F^j_i:\mathcal{B}P_i\rightarrow \mathcal{L}\}^{1\leq j\leq 2k(i)+1}_{1\leq i\leq n}$ in the following way. For $j=1,\cdots , k(i)$ $F^{2j-1}_i:\mathcal{B}P_i\rightarrow\mathcal{L}$, $x\mapsto \widehat{\psi ^{j-1}_i}\circ\psi ^{j-2}_i\circ\cdots\circ\psi ^1_i(x)$, $\bullet\mapsto S$, $F^{2j}_i:\mathcal{B}P_i\rightarrow\mathcal{L}$, $x\mapsto \widehat{\psi ^{j-1}_i}\circ\psi ^{j-2}_i\circ\cdots\circ\psi ^1_i(x)$, $\bullet\mapsto P_i^j$,$F^{2j+1}_i:\mathcal{B}P_i\rightarrow\mathcal{L}$, $x\mapsto \widehat{\psi ^{j}_i}\circ\psi ^{j-1}_i\circ\cdots\circ\psi ^1_i(x)$, $\bullet\mapsto S$. It follows from the existence of the linking system that we can find lifts of the inclusion $\{\iota _{P_i^j,S}\}$ and lifts of the morphisms $\{\widehat{\psi ^j_i}\}$ which will be denoted $\{\widetilde{\psi ^j_i}\}$ such that for all $i=1,\cdots n$ and for all $l=1,\cdots 2k(i)+1$ the functors $F_i^{2j}$ and $F_i^{2j+1}$ commute via the following natural transformations respectively. The functor $F_i^{2j}$ commutes to $F^{2j-1}_i$ via $\bullet\mapsto\iota _{P_i^j}$ and $F^{2j}_i$ commutes to $F^{2j+1}_i$ via $\bullet\mapsto\widetilde{\psi ^j_i}$. Note that we have an algebraic map $q:H^*(B\mathcal{G})\rightarrow H^*(|\mathcal{L}|)$ as before.
The induced diagram on classifying spaces 
commutes up to homotopy for all $i=1,...,n, j=1,2$. We have $s:B\mathcal{G}= \underset{\mathcal{D}}{hocolim^{(1)}}(F)\rightarrow |\mathcal{L}|$
which induces $s^*:H^*(|\mathcal{L}|;\mathbb{F}_p)\longrightarrow H^*(B\mathcal{G};\mathbb{F}_p)$. Analogously to the case of Robinson models we have the following commutative diagram
\xymatrix@R=7pt@C=7pt{ {BS}\ar[0,2]^{B(\delta _S)}\ar[1,1]_{Bincl}&{}&{|\mathcal{L}|}\\
&{B\mathcal{G}}\ar[-1,1]_{s}&\\} inducing a commutative diagram of unstable algebras over the Steenrod algebra
\xymatrix@R=7pt@C=7pt{ {H^*(BS)}&{}&{H^*(|\mathcal{L}|)}\ar[0,-2]_{B(\delta _S)^*}\ar[1,-1]^{s^*}\\
&{H^*(B\mathcal{G})}\ar[-1,-1]^{q}&\\} showing that $q\circ s^*=id_{H^*(|\mathcal{L}|)}$.$\Box$\\

\begin{Theorem}
Let $(S,\mathcal{F},\mathcal{L})$ be a $p-$local finite group and $\mathcal{G}$ a model of Leary-Stancu type for $\mathcal{F}$. Then
$|\mathcal{L}|\pcom$ is a stable retract of $B\mathcal{G}\pcom$.
\end{Theorem}
\underline{Proof:} 
The diagram
\xymatrix@R=9pt@C=9pt{&{\Sigma ^{\infty}BS\pcom}\ar[1,1]^{\Sigma ^{\infty}Bincl\pcom}\ar[1,-1]_{\Sigma ^{\infty}B(\delta _S)\pcom}&{}\\
{\Sigma ^{\infty}|\mathcal{L}|\pcom}&&{\Sigma ^{\infty}B\mathcal{G}\pcom}\ar[0,-2]^{\Sigma ^{\infty}r\pcom}\\
}
commutes where $s$ is the map constructed in the proof of Theorem $4.8$.
Since $\Sigma ^{\infty}B(\delta _S)\pcom\circ \sigma _{\mathcal{F}} =\Sigma ^{\infty}s\pcom\circ\Sigma ^{\infty}Bincl\pcom\circ \sigma _{\mathcal{F}}$ we have $|\mathcal{L}|\pcom$ is a stable retract of $B\mathcal{G}\pcom$.$\Box$\\
\begin{Theorem}
Let $\mathcal{F}$ be a saturated fusion system over the finite $p-$group $S$ and $\mathcal{G}$ a model of Leary-Stancu type for $\mathcal{F}$ with set of automorphisms $\Phi$. Then $H^*(B\mathcal{G})$ is noetherian if and only if $\Phi$ is finite.
\end{Theorem}
\underline{Proof:} Note that we have a map $B\mathcal{G}=\underset{\mathcal{D}}{hocolim(F)}\rightarrow |\mathcal{L}|$ where $F$ and $\mathcal{D}$ are as defined in the proof of Theorem 4.8. for the model of Leary-Stancu type $\mathcal{G}$. Again $H^*(B(P_i))$ is finitely generated over $H^*(|\mathcal{L}|)$ for all $i=1,...,n$, and $H^*(|\mathcal{L}|)$ is noetherian. Therefore the Bousfield-Kan spectral sequence for $H^*(B\mathcal{G})$ is a spectral sequence of finitely generated $H^*(|\mathcal{L}|)-$modules, the $E_2$ term with $E_2^{s,t}= \underset{\mathcal{C}}{lim^s}H^t(F(-);\mathbb{F}_p)$ is concentrated in the first two columns so $E_{2}=E_{\infty}$ for placement reasons. Therefore $H^*(B\mathcal{G})$ is a finitely generated module over $H^*(|\mathcal{L}|)$ and in particular noetherian if and only if $\Phi$ is finite. $\Box$\\[0.1cm]

\bibliographystyle{amsplain}
Laboratoire Analyse, Géométrie et Applications\\
UMR 7539\\
Institut Galilée\\
Université Paris 13\\
99 avenue J.B. Clément\\
93430 Villetaneuse\\
FRANCE \\
\textit{email:} seeliger@math.univ-paris13.fr
\end{document}